\documentclass[conference]{IEEEtran}
\IEEEoverridecommandlockouts
\usepackage{cite}
\usepackage{amsmath,amssymb,amsfonts}
\usepackage[ruled,vlined]{algorithm2e}
\usepackage{graphicx}
\usepackage{textcomp}
\usepackage{xcolor}
\usepackage[utf8]{inputenc}
\usepackage{url}

\begin{document}

\title{%
Image edges resolved well when using an overcomplete piecewise-polynomial model
\thanks{Research described in this paper was financed by the National Sustainability Program under grant LO1401 and by the Czech Science Foundation under grant no.\,GA16-13830S.
For the research, infrastructure of the SIX Center was used.}
}

\author{
\IEEEauthorblockN{Michaela Novosadová}
\IEEEauthorblockA{\textit{Signal Processing Laboratory} \\
\textit{Brno University of Technology}\\
Brno, Czech Republic \\
novosadova@phd.feec.vutbr.cz}
\and
\IEEEauthorblockN{Pavel Rajmic}
\IEEEauthorblockA{\textit{Signal Processing Laboratory} \\
\textit{Brno University of Technology}\\
Brno, Czech Republic \\
rajmic@feec.vutbr.cz}
}

\maketitle

\begin{abstract}
Used in the paper is an overcomplete piecewise-polynomial image model incorporating sparsity.
The paper shows that using such a model, the edges in the image can be resolved robustly with respect to noise.
Two variants of the proposed approach are both shown to be superior to the use of the classic edge detecting kernels.
The proposed method is in turn also suitable for image segmentation.
%
\end{abstract}


\newcommand{\todo}[1]{\textcolor{red}{#1}}
\newcommand\Rset{\mathbb{R}}
\newcommand{\supp}{\mathop{\operatorname{supp}}}
\newcommand{\vectorize}{\mathop{\operatorname{vec}}}
\newcommand{\reshape}{\mathop{\operatorname{reshape}}}
\newcommand{\soft}{\mathop{\operatorname{soft}}}
\newcommand{\prox}{\mathop{\operatorname{prox}}}
\newcommand{\proj}{\mathop{\operatorname{proj}}}
\newcommand{\diag}{\mathop{\operatorname{diag}}}
\newcommand{\argmin}{\mathop{\operatorname{arg~min}}}
\newcommand{\TV}{\mathop{\operatorname{TV}}}
\def\x{\vect{x}}
\def\q{\vect{q}}
\def\p{\vect{p}}
\def\e{\vect{e}}
\def\u{\vect{u}}
\def\y{\vect{y}}
\def\z{\vect{z}}
\def\f{\vect{f}}
\def\d{\vect{d}}
\def\w{\vect{w}}
\def\A{\matr{A}}
\def\U{\matr{U}}
\def\V{\matr{V}}
\def\S{\matr{S}}
\def\X{\matr{X}}
\def\B{\matr{B}}
\def\N{\matr{N}}
\def\O{\matr{O}}
\def\R{\matr{R}}
\def\P{\matr{P}}
\def\I{\matr{I}}
\def\D{\matr{D}}
\def\W{\matr{W}}
\def\Z{\matr{Z}}
\def\L{\matr{L}}
\newcommand{\norm}[1]{\|#1\|}
\newcommand{\Norm}[1]{\left\|#1\right\|}
\newcommand{\abs}[1]{\left\vert#1\right\vert}
\newcommand{\transp}[1]{#1^\top}
\newcommand{\vect}[1]{\mathbf{#1}} 
\newcommand{\matr}[1]{\mathbf{#1}} 
\newcommand{\qm}[1]{``#1''}  

\section{Introduction}
Segmentation of signals, and images in particular, is an important
field of research with a huge number of applications.
The individual segments of an image are distinguishable due to the fact
that they typically exhibit
(in a~wide sense)
a~consistent character,
and this character changes at the border of segments.
Models beyond the individual segmentation algorithms most often rely on the mathematical modeling
of what is understood under the vague term \qm{character};
however, there are also approaches based on machine-learning segmentation from training data 
\cite{Milletari2016:Fully.convolutional.NN.Volumetric.segmentation,%
Fritscher2016:.DNN.Fast.Segmentation.3D.medical}.

Within the approaches based on a~priori modeling, one can distinguish explicit and implicit types of models.
In the \qm{explicit} type, the signal is modeled such that it is a composition of sub-signals which can often be expressed analytically
\cite{GiryesEladBruckstein2015:overcomplete.sparse.variational,%
Shem-Tov2012:globally.optimal.local.modeling,%
RajmicNovosadDankova2017:Piecewise-polyn.signal.segment_Kybernetika,%
OngieJacob2015:Recovery.discontinuous.signals.group.sparse,%
SelesnickArnoldDantham2012:polynomial.step.discontin,%
ZhangGengLai2015:multiple.change.sparse.group.lasso}.
In the \qm{implicit} type,
the signal is characterized by features that are derived from the signal by using a~suitable linear operator
\cite{BleakleyVert2011:group.fused.Lasso,%
KimKohBouydGorinevsky2009:l1.trend.filtering,%
Selesnick2017:Sparsity.assisted.signal.smoothing.revisited,
Tibshirani2014:Adaptive.piecewise.polynomial%
}.
The above differences are analogous to the so-called \qm{synthesis} and \qm{analysis} approaches, respectively,
recognized in the sparse signal processing literature
\cite{Elad05analysisversus,NamDaviesEladGribon2013:CosparseAnalysisModel}.
The two types of models are significantly different in the way the signal is treated;
however, connections between them can be found, for example see
\cite{UnserFageotWard2017:Splines.universal.solutions.GTV}.

In the paper, we accept the assumption that images consist of piecewise-polynomial patches,
while it is typically sufficient to consider quadratic polynomials 
\cite{Elad2005:Cartoon.and.texture-MCA,BrediesHoller2015:TGV-based.framework.partI,HollerKunisch2014:Infimal.convolution.TV.type.video.image,Knoll2011:second.order.tv,KutyniokLim2011:Shearlets.are.optimal}.
In that regard, the segment border is a~place in the image where the polynomial characterization is subjected to a~change,
while within the segments, the respective representation stays steady.

With noise preset in the image, the segmentation success naturally decreases.
In the paper, we propose two methods based on a~single model and show that both of them perform better in edge detection than standard methods do.
Our work was inspired by
\cite{GiryesEladBruckstein2015:overcomplete.sparse.variational},
where the authors use a similar model, which they approximately solve using a~greedy algorithm. 
Contrariwise, our approach relies on the convex relaxation of the sparsity measure
\cite{DonohoElad2003:Optimally}.
As far as we know, the present paper is the first to use convex relaxation for modeling edges in images in connection with the piecewise-polynomial assumption. 

In our recent works, we dealt with various segmentation algorithms for 1D signals, all of them following the piecewise-polynomial model
\cite{NovosadovaRajmic2016:Piecewise-polynom.curve.fitting,
RajmicNovosadDankova2017:Piecewise-polyn.signal.segment_Kybernetika,
NovosadovaRajmicSorel2018:Orthogonal.segmentation}.
In the present work, we generalize the 1D model to 2D, derive a new numerical solver and adapt the way the model output is processed.
Note that in terms of the above categorization, the proposed model belongs to the \qm{explicit} class.


\section{Overcomplete piecewise-polynomial model}
\label{Sec:Model}
The images are modeled as being composed of non-overlapping polynomial patches.
To better understand the way the two-dimensional polynomials are constructed,
we first present the one-dimensional model and then generalize it.
For more details on the 1D model, we refer to
\cite{NovosadovaRajmicSorel2018:Orthogonal.segmentation}.

\subsection{One-dimensional model}
\label{Sec:1D.model}
In continuous time, a~polynomial of degree $K$ is representable as a~linear combination of basis polynomials:
\begin{equation}
	y(t) = x_0 p_0(t) + x_1 p_1 (t) + \ldots + x_K p_K(t) ,\quad  t\in\Rset,
\end{equation}
where $x_k$ are the expansion coefficients in such a~basis.
The system $\{p_k(t)\}_{k=0}^K$ is only required to form a basis,
but additional properties such as orthogonality can sometimes be beneficial,
see \cite{NovosadovaRajmicSorel2018:Orthogonal.segmentation}.
%
In a~discrete setting, the individual elements of a~polynomial signal
are represented as
\begin{equation}
	\label{Eq:Discrete.polynomial.signal}
	y[n] = x_0 p_0[n] + x_1 p_1[n] + \ldots + x_K p_K[n],\quad n=1,\ldots,N.
\end{equation}

Given $\{p_k(t)\}_{k=0}^K$,
every signal given by \eqref{Eq:Discrete.polynomial.signal} is uniquely determined by the~set of coefficients $\{x_k\}_{k=0}^K$.
The uniqueness of such an expansion is broken when we introduce a~time index also in these coefficients, allowing them to change in time:
\begin{equation}
	y[n] = x_0[n] p_0[n] + x_1[n] p_1[n] + \ldots + x_K[n] p_K[n].
\end{equation}
This makes the representation of the signal $\{y[n]\}_{n=1}^N$ extremely overcomplete,
nevertheless, this way of handling plays a~key role in the model below.
It will be convenient to write the latest relation in a~more compact form, for which we introduce the notation
\begin{equation}
	\y = \begin{bmatrix} y[1]\\ \vdots\\ y[N] \end{bmatrix}\!,\
	\x_k = \begin{bmatrix} x_k[1]\\ \vdots\\ x_k[N] \end{bmatrix}\!,\
	\p_k = \begin{bmatrix} p_k[1]\\ \vdots\\ p_k[N] \end{bmatrix}
\end{equation}
for $k=0,\ldots, K$.
After this, we can write
\begin{equation}
	\label{eq:1D.overcomplete.polynomial.model.as.elementwise.multiplic}
	\y = \p_0\odot\x_0 + \ldots + \p_K\odot\x_K
\end{equation}
with $\odot$ denoting the elementwise (Hadamard) product.
To convert \eqref{eq:1D.overcomplete.polynomial.model.as.elementwise.multiplic}
to the form involving the standard matrix--vector product, we define
\begin{equation}
	\P_k = \diag (\p_k) = \begin{bmatrix} p_k[1] & & 0\\  & \ddots & \\  0 &  & p_k[N] \end{bmatrix}\!,\
	k=0,\ldots, K,
\end{equation}
allowing us to write 
\begin{equation}
	\y = \P_0 \x_0 + \ldots + \P_K \x_K
\end{equation}
or, even more compactly,
\begin{equation}
	\label{eq:Polynomial.model.in.short.form}
	\y = \P \x = \left[\P_0 | \cdots | \P_K\right]
	\begin{bmatrix}\x_0\\[-1ex] \rotatebox{90}{$|$} \\[-1ex] \vdots \\[-1ex] \rotatebox{90}{$|$} \\[-1ex] \x_K\end{bmatrix}.
\end{equation}

Such a description of a~signal of dimension $N$ is superfluous in terms of the number of parameters;
when the polynomials are fixed (which is assumed in the rest of the paper),
the total number of $(K+1)N$ parameters are used to characterize the signal.
Nevertheless, assume now that $\y$ is a~\emph{piecewise} polynomial and that it consists of $S$ independent segments.
\emph{Each} segment $s\in\{1,\ldots,S\}$ is then described by $K+1$ polynomials.
In the above context, this can be achieved by letting \emph{all} the vectors $\x_k$ be constant within the particular segments, with a~change at the segment border.
%
%
Note that the positions of the segment change-points are unknown.

The described model motivates the following idea.
The finite difference operator $\nabla$ applied to the vectors $\x_k$
produces vectors that are sparse---actually, $\nabla\x_k$ has $S-1$ nonzeros at maximum.
Moreover, the nonzero components of all $\nabla\x_k$ occupy the same positions across $k=0,\ldots, K$.
This property justifies the use of the group-sparse model introduced later.

\subsection{Two-dimensional model}
\label{Sec:2D.model}
The natural extension to 2D is done through the Kronecker product of 1D vectors.
For brevity, the derivation of the polynomial model will be introduced 
under the (natural) assumption that the degrees of polynomials in the vertical and horizontal directions are identical.

\begin{figure}
	\includegraphics[width=.48\columnwidth]{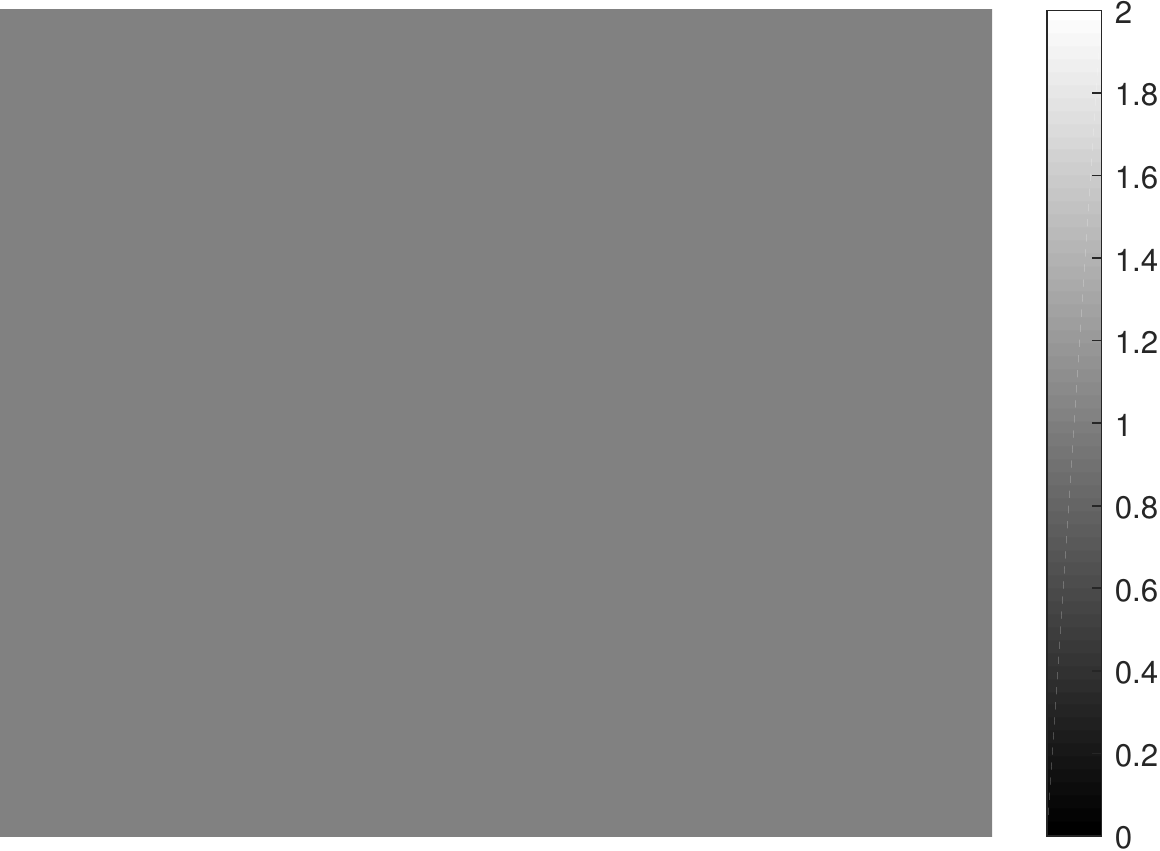}%
	\hfill%
	\includegraphics[width=.48\columnwidth]{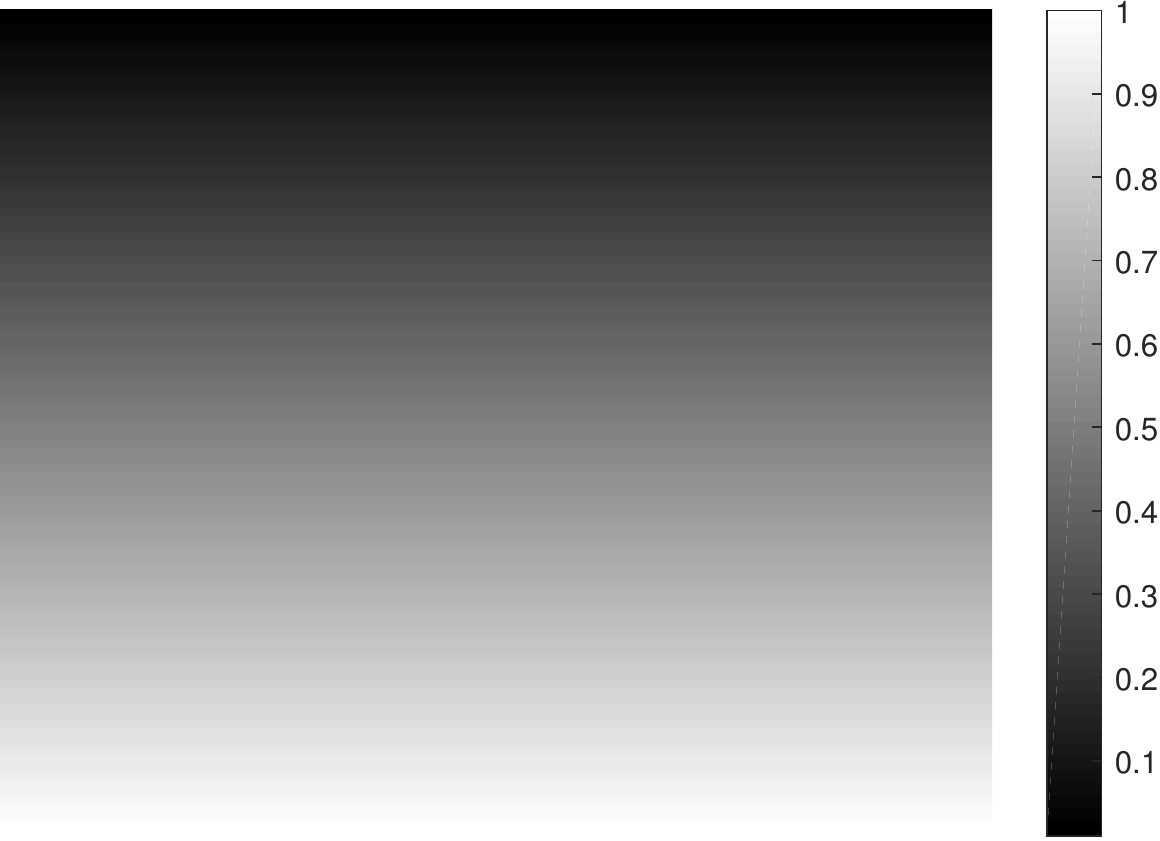}\\[2pt]
	\includegraphics[width=.48\columnwidth]{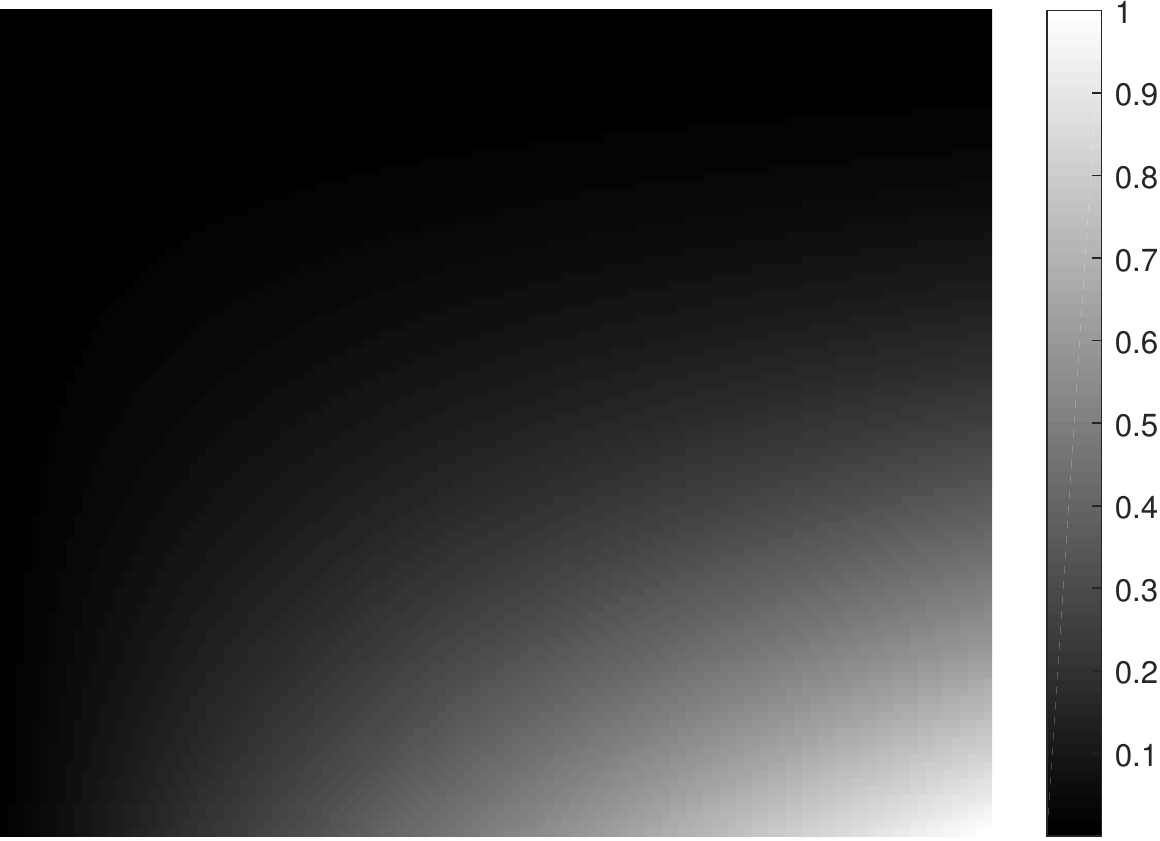}%
	\hfill%
	\includegraphics[width=.48\columnwidth]{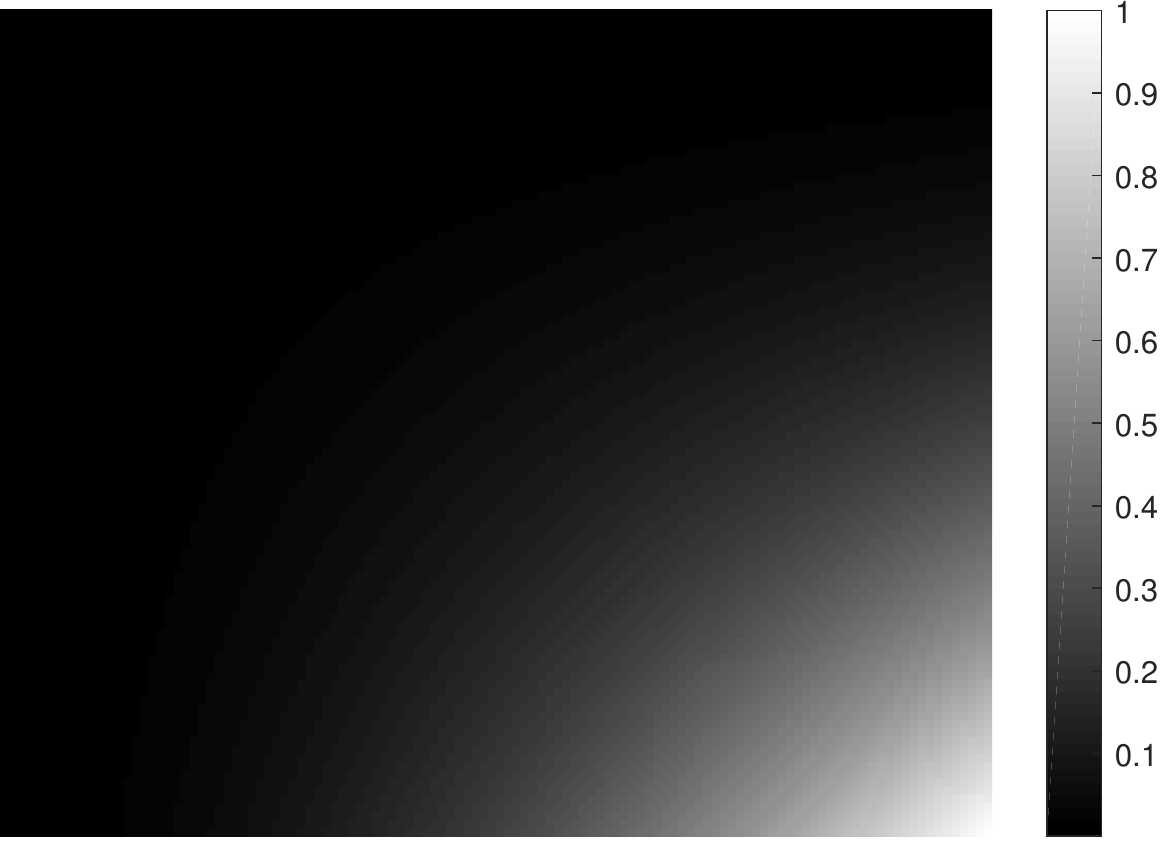}
	\caption{Polynomials $\p_{00}$, $\p_{10}$, $\p_{21}$, $\p_{22}$ of size $100\times 120$, generated from the standard basis of degree $K=2$
	(i.e.\ from $1,t,t^2$).
	}
	\label{fig:Examples.of.2D.polynomials}
\end{figure}

The 2D polynomials are formed by all the combinations
\begin{align*}
	\p_{00} &= \p_{0\mathrm{v}} \cdot \transp{\p_{0\mathrm{h}}} \\
	\p_{01} &= \p_{0\mathrm{v}} \cdot \transp{\p_{1\mathrm{h}}} \\
				&\ \vdots \\
	\p_{KK} &= \p_{K\mathrm{v}} \cdot \transp{\p_{K\mathrm{h}}}
\end{align*}
making altogether $(K+1)^2$ polynomial basis images.
The subscripts \qm{v} and \qm{h} distinguish the generating 1D polynomials
in the vertical and horizontal direction, respectively.
See Fig.\,\ref{fig:Examples.of.2D.polynomials} for examples.
It is easy to show that when both the sets of vertical and horizontal prototypes form two respective orthonormal bases,
then the new 2D system, too, is an orthonormal basis.

Let the image signal be $\y\in\Rset^{M\times N}$ with entries
$y[m,n]$, $m=1,\ldots,M$, $n=1,\ldots,N$.
As in \eqref{eq:1D.overcomplete.polynomial.model.as.elementwise.multiplic},
the image is modeled as
\begin{equation}
	\label{eq:2D.overcomplete.polynomial.model.as.elementwise.multiplic}
	\y = \p_{00}\odot\x_{00} + \p_{01}\odot\x_{01} + \ldots + \p_{KK}\odot\x_{KK},
\end{equation}
where $\p_{k\ell}$ are basis images of size $M\times N$ as defined above, and
$\x_{k\ell}$ are matrices of parameters, all of the same size $M\times N$.

We need to identify an image with its vectorized form;
this way of perception is necessary for us to be able again to express the model as a~matrix--vector multiplication.
Let the vectorization operator $\vectorize(\cdot)$ stack image columns, one after another, into a single column.
The model \eqref{eq:2D.overcomplete.polynomial.model.as.elementwise.multiplic}
reads
\begin{equation}
	\vectorize(\y) = \P_{00} \cdot \vectorize(\x_{00}) + \ldots + \P_{KK} \cdot \vectorize(\x_{KK})
\end{equation}
when we define $\P_{k\ell}=\diag (\vectorize (\p_{k\ell}))$.
In short, convenient form,%
\begin{equation}%
	\label{eq:Polynomial.model.in.short.form}
	\vectorize(\y) = \P \x = \left[\P_{00} | \cdots | \P_{KK}\right]
	\begin{bmatrix}\vectorize(\x_{00})\\[-1ex] \rotatebox{90}{$|$} \\[-1ex] \vdots \\[-1ex] \rotatebox{90}{$|$} \\[-.5ex] \vectorize(\x_{KK})\end{bmatrix}.
\end{equation}

\subsection{Basic idea}
Following the just described image-generation model relying on $\P$,
our method is based on the convex optimization problem presented in the next section.
After the numerical solver has found an optimal solution to the problem,
the solution is used in two possible ways to resolve the edges of the image.
Both variants include the thresholding operator as the final step.

\subsection{Optimization problem}
As in the above case, in the 2D setting, too, the parameterization coefficients $\{\x_{k\ell}\}_{k,\ell}$
can be chosen such that they describe the signal perfectly and, at the same time, they are sparse under the difference operator
$\nabla$.
However, the difference operator has two forms now---the horizontal and the vertical one.

This, together with the assumption that the observed signal is corrupted by an i.i.d.\ Gaussian noise,
motivates the following formulation of the denoising problem:
%
\begin{multline}
	\label{Eq:Problem.formulation}
	\hat{\x}
	=
	\argmin_{\x}
	\left\{ \norm{\L_\mathrm{v}\x}_{21} + \lambda\norm{\L_\mathrm{h}\reshape(\x)}_{21} \right\} \\
	\text{ s.t.\ } \norm{\!\vectorize(\y)-\P\x}_2 \leq \delta.
\end{multline}
%
In this optimization program,
the operators $\L_\mathrm{v}$ and $\L_\mathrm{h}$ represent,  respectively, the stacked vertical and horizontal differences,
such that
%
\begin {align}
	\label{eq:Definition.of.linear.operator}
	\L_\mathrm{v} &=
		\begin{bmatrix}
			\ \nabla & \cdots   & 0 \\
												& \ddots   &\\
			0                 & \cdots   & \nabla
		\end{bmatrix},
\end {align}
with $\nabla: \Rset^{M}\mapsto\Rset^{M-1}$ being the standard columnwise difference operator.
This operator is repeated $(K+1)^2$ times---the number of parameter maps---in \eqref{eq:Definition.of.linear.operator}.
The operator $\L_\mathrm{h}$ is defined identically as $\L_\mathrm{v}$ in \eqref{eq:Definition.of.linear.operator},   but with the difference that $\nabla$ is now the mapping $\nabla: \Rset^{N}\mapsto\Rset^{N-1}$.
The operator $\reshape(\cdot)$ is used to emphasize that it is necessary to tweak the entries of $\x$ such that the horizontal differences are computed using the linear operator $\L_\mathrm{h}$.

The first term of \eqref{Eq:Problem.formulation} is the penalty.
Piecewise constant $\x_{k\ell}$ 
suggest that these vectors are sparse under the difference operation. 
The usual convex substitute of the true sparsity measure is the $\ell_1$-norm 
\cite{CandesWakinBoyd2008:enhanced.reweighted.l1,DonohoElad2003:Optimally}.
Recall, however, that in addition to sparsity, the nonzero differences should appear at the same spots,
as argumented in Sec.~\ref{Sec:1D.model}.
Therefore, our optimization should contain a~penalty that promotes \emph{group} sparsity, whose convex substitute is the $\ell_{21}$-norm
\cite{Kowalski.Torresani2009:Structured.sparsity.from.mixed}.
The $\ell_{21}$-norm acts on a~matrix 
and is formally defined by%
\begin{align}%
	\label{eq:Definition.of.ell21-norm}
	\norm{\Z}_{21} &=
	\norm{\Z_{1,:}}_2 + \ldots + \norm{\Z_{p,:}}_2 \,,
\end{align}
i.e.\ the $\ell_2$-norm is applied to the individual rows of $\Z$ and the resulting vector is measured by the $\ell_1$-norm.
In our case, the matrix is formed such that it is of size $(M-1)(N-1)\times (K+1)^2$ and as such, the $\ell_{21}$ penalty promotes sparsity across differentiated pixels, in line with the justification of the model.
For details, see 
\cite{RajmicNovosadDankova2017:Piecewise-polyn.signal.segment_Kybernetika}.

The scalar parameter $\lambda$ optionally weights the vertical vs.\ horizontal differences.
This parameter can differ from unity in the case that the user wants to emphasize horizontal over vertical structures
and vice versa. 
Note, however, that this parameter should not change when the size of the image should \qm{grow} proportionally in both directions, since in such a case both the penalties grow proportionally as well;
only $\delta$ should scale with the image size.

The second term in \eqref{Eq:Problem.formulation} is the data fidelity term.
The Euclidean
norm
reflects the fact that gaussianity of noise is assumed. 
Consider however, that $\delta$ accounts not only for the noise power, but also for model imperfections,
since real images do not follow the piecewise-polynomial assumption exactly.
As such, this parameter has to be tuned.
Finally, the vector $\hat{\x}$ contains the obtained optimizers.

\section{Numerical solver}
The optimization problem 
\eqref{Eq:Problem.formulation}
requires minimizing a sum of convex functions within a~convex set of feasible solutions.
Proximal splitting algorithms are suitable for such a~class of tasks.
We utilize the Condat proximal algorithm
\cite{Condat2014:Generic.proximal.algorithm}, 
which is able to solve (among more general formulations) problems of the type
\begin{equation}
	\label{Eq:minimization.problem.2}
	\text{minimize}\ h_1(\L_1\x) + h_2(\L_2\x) + h_3(\L_3\x),
\end{equation}
over $\x$, where $h_i$ are convex functions and $\L_i$ are linear operators. 

The connection between \eqref{Eq:minimization.problem.2} and \eqref{Eq:Problem.formulation}
is provided via assigning
$h_1 = \norm{\cdot}_{21}$,
$h_2 = \lambda \norm{\cdot}_{21}$,
$\L_1 = \L_\mathrm{v}$,
$\L_2 = \L_\mathrm{h}\reshape(\cdot)$.
To fit the formulation \eqref{Eq:minimization.problem.2}, the feasible set is recast in the 
unconstrained form using the indicator function
$h_3 = \iota_{\{\z:\, \norm{\y-\z}_2 \leq \delta\}}$,
where $\iota_C$ denotes the indicator function of a~convex set $C$ 
\cite{combettes2011proximal}.
Finally, $\L_3 = \P$.
 
The algorithm that solves
\eqref{Eq:Problem.formulation}
is described in Alg.\,\ref{fig:Condat.algorithm}.
Therein, two operators are involved which have not been introduced yet:
The operator $\soft^\text{group}_{\tau}$ performs \qm{group} soft thresholding with the threshold $\tau$,
where there are $(M-1)(N-1)$ disjoint groups, each containing $(K+1)^2$ coefficients.
The projector $\proj_{B_2(\y,\delta)}$ finds the closest point in the $\ell_2$-ball
$\{\z: \norm{\y-\z}_2 \leq \delta\}$ with respect to the input point.
Both the operations are computationally cheap, for details see for example
\cite{RajmicNovosadDankova2017:Piecewise-polyn.signal.segment_Kybernetika}.

\begin{algorithm}[t]
\SetAlgoVlined
\label{fig:Condat.algorithm}	
	
	\KwIn{%
		$\P$, 
		$\y$,
		$\delta$} 
	\KwOut{$\hat{\x}=\x^{(i+1)}$}
	
			Set the parameters $\xi, \sigma >0$ and $\rho \in (0, 2)$. 
		
		%
			Set the initial primal variable $\x^{(0)}$ 
			and dual variables $\u_1^{(0)},\u_2^{(0)},\u_3^{(0)}$. 
		

	\For{$i=0,1,\dots$}{
$\bar{\x}^{(i+1)} = \x^{(i)}-\xi
		\left(
			\transp{\L_\mathrm{v}} \u_1^{(i)} + \transp{\reshape}\transp{\L_\mathrm{h}} \u_2^{(i)} + \transp{\P}\u_3^{(i)}
		\right)$

$\x^{(i+1)} = \rho\,\bar{\x}^{(i+1)} + (1-\rho)\x^{(i)}$

$ \p_1 = \u_1^{(i)} + \sigma\, \L_\mathrm{v} (2\bar{\x}^{(i+1)}-\x^{(i)})$

$\bar{\u}_1^{(i+1)} =	\p_1 - \soft^\text{group}_{1}(\p_1)$

$\u_1^{(i+1)}=\rho\,\bar{\u}_1^{(i+1)}+(1-\rho)\u_1^{(i)}$

$ \p_2 = \u_2^{(i)} + \sigma\, \L_\mathrm{h} \reshape (2\bar{\x}^{(i+1)}-\x^{(i)})$

$\bar{\u}_2^{(i+1)} =	\p_2 - \soft^\text{group}_{\lambda}(\p_2 )$

$\u_2^{(i+1)}=\rho\,\bar{\u}_2^{(i+1)}+(1-\rho)\u_2^{(i)}$

$ \p_3 = \u_3^{(i)} + \sigma\, \P(2\bar{\x}^{(i+1)}-\x^{(i)})$

$\bar{\u}_3^{(i+1)}= \p_3 - \sigma \,\proj_{B_2(\y,\delta)}(\p_3/\sigma)$

$\u_3^{(i+1)}=\rho\,\bar{\u}_3^{(i+1)}+(1-\rho)\u_3^{(i)}$

	}
	\KwRet{${\x}^{(i+1)}$}
	
	\caption{The Condat Algorithm solving \eqref{Eq:Problem.formulation}} 
	
\end{algorithm}

The convergence  of the algorithm is guaranteed when it holds
$\xi\sigma \norm{\transp{\L_1}\! \L_1+\transp{\L_2}\! \L_2+\transp{\L_3}\! \L_3}\leq 1$,
where $\norm{\cdot}$ denotes the operator norm.
The operator norm of the sum can be bounded such that
$
\norm{\transp{\L_1}\! \L_1+\transp{\L_2}\! \L_2+\transp{\L_3}\! \L_3}
\leq \norm{\L_1}^2 + \norm{\L_2}^2 + \norm{\L_3}^2
$.
We have 
\begin{align*}
\begin{split}
	\norm{\L_1}^2 &= \norm{\L_\mathrm{v}}^2
			 = \max_{\norm{\x}_2=1} \norm{\L_\mathrm{v}\x}^2_2 \\
			 &=	\max_{\norm{\x}_2=1} \left( \sum_{k=0}^K \sum_{\ell=0}^K \Norm{\nabla\vectorize(\x_{k\ell})}^2_2 \right) \\
			 &\leq \sum_{k=0}^K \sum_{\ell=0}^K \left( \max_{\norm{\x}_2=1}  \Norm{\nabla\vectorize(\x_{k\ell})}^2_2\right) \\
			 &\leq \sum_{k=0}^K \sum_{\ell=0}^K \, \norm{\nabla}^2
			 = \norm{\nabla}^2 (K+1)^2 \,\leq\, 4 (K+1)^2
\end{split}
\end{align*}
%
and the same holds for $\norm{\L_2}$.
As regards $\norm{\L_3}=\norm{\P}$,
we have
$\norm{\P}^2 = \norm{\P\transp{\P}\!}$
and thus it suffices to find the maximum eigenvalue of
$\P\transp{\P}\!$.
Since $\P$ has a~multi-diagonal structure (cf.\ \eqref{eq:Polynomial.model.in.short.form}),
$\P\transp{\P}\!$ is diagonal and thus it suffices to find the maximum on its diagonal.
Altogether, the convergence is guaranteed when
$\xi\sigma \left( \max\diag(\P\transp{\P}\!) + 8(K+1)^2 \right) \leq 1$.


\section{Experiments}
The goal of the experiment was to establish whether the proposed
problem
\eqref{Eq:Problem.formulation} can help in indicating edges in noisy images better than the standard edge detection methods.

The experiment was designed similar to that of the evaluation in
\cite{GiryesEladBruckstein2015:overcomplete.sparse.variational}.
The starting signals for the test were the standard grayscale photographs taken from Matlab, the range of whose values was 0 to 255.
The images were contaminated with zero mean i.i.d.\ gaussian noise.

\subsubsection{Coins, lower noise}
The first experiment was performed with the photo depicting a~number of coins,
see Fig.\,\ref{fig:Coins.processed}.
The image is of size 246$\times$300\,px.
We superposed the noise over it with the standard deviation $\sigma=20$.

As the ground truth edges we consider the result of the Sobel operator
\cite{Burger.Burge_2010:Digital.image.processing}
applied to the original, clean image, subsequent normalization of the result such that the maximum absolute value is one, and finally the  thresholding of the absolute values.
The output, i.e.\ the binary map, is supposed to indicate the true positions of edges.
We chose a~threshold value of 0.15 for this particular image, which retained the important details at the sides of the coins.
See Fig.\,\ref{fig:Coins.processed}.

As the standard way of edge detection we consider the application of the Sobel operator to the noisy image and thresholding
(the normalized absolute values).
In this process, we tuned the threshold to a compromise such that false edges do not appear too often
while some of the details are still preserved in the coins.
In our particular case, the threshold was 0.24.

The newly proposed approach was as follows:
In the optimization problem \eqref{Eq:Problem.formulation}, the standard basis of degree $K=2$ was used as the generating polynomial basis,
corresponding to the examples in Fig.\,\ref{fig:Examples.of.2D.polynomials}.
The basis images 
$\p_{00},\ldots,\p_{22}$ were used to form the model matrix $\P$ as described in Sec.\ \ref{Sec:2D.model}.
The model error $\delta=4\,000$ was allowed.

We let Alg.\,\ref{fig:Condat.algorithm} run for 500 iterations with this setting
(which took 3 minutes on a~PC equipped with Intel i7 processor and 16\,GB of RAM).
This way we obtained the solution to \eqref{Eq:Problem.formulation}, namely $\hat\x$. 
This vector is used in two distinct ways to get an estimate of the image edges.

The first way simply synthesizes the denoised image, i.e.\ it forms the image by
$\hat{\y} = \P \hat\x$ and reshapes the resulting column into a matrix of the size 264$\times$300.
Then, the above described convolution with the Sobel kernel and thresholding is applied as in the above cases.
Again, the threshold was found (0.14) to produce the best-looking compromise between false alarms and edge preservation.

The second way exploits the multiple information contained in $\hat\x$, in contrast to using the single synthesized image alone.
In particular, recall that the edges are (at least in the model) indicated in all the sets of parameters $\x_{k\ell}$ in the same places,
cf.\ \eqref{eq:2D.overcomplete.polynomial.model.as.elementwise.multiplic}.
The space of $\x$ is much larger than the image space and therefore the parameters could potentially be helpful in detecting the edges.
We first compute the Sobel operator on all the individual parameter maps $\hat\x_{k\ell}$.
Then we compute pixelwise $\ell_2$-norms along the $(K+1)^2$ parameters.
The output is a single matrix and is subject to thresholding (in our case with the threshold 0.12).

In total, we compare three methods (Sobel and that proposed in two versions) with the \qm{ground truth}.
The results can be seen in Fig.\,\ref{fig:Coins.processed}.

\subsubsection{Coins, higher noise}
The same experiment but with greater noise is presented in Fig.\,\ref{fig:Coins.processed.more.noise}
(the respective thresholds are 0.15, 0.29, 0.16 and 0.18; $\delta=6\,000$).

\subsubsection{Coins, higher noise, orthonormal basis}
\label{sec:coins.ortho}
The same image was processed using the orthonormal basis instead of the standard basis
\cite{NovosadovaRajmicSorel2018:Orthogonal.segmentation}.
This time, the model error was set to $\delta=8\,000$ and the respective thresholds are 
0.15, 0.29, 0.14 and 0.14.
To reach convergence, we had to let Alg.\,\ref{fig:Condat.algorithm} run for 10\,000 iterations.
The experiment is presented in Fig.\,\ref{fig:Coins.processed.more.noise.orthonormal}.

\subsubsection{Cameraman}
The cameraman image of the size 256\,$\times$\,256 (Fig.\,\ref{fig:Cameraman.processed}) was processed using the standard basis and with
the model error $\delta=7\,000$ and the respective thresholds are 
0.25, 0.29, 0.18 and 0.16.
To reach convergence, we had to let Alg.\,\ref{fig:Condat.algorithm} run for 500 iterations.

\subsubsection{Rice}
The rice image of the size 256\,$\times$\,256 (Fig.\,\ref{fig:Rice.processed}) was processed using the standard basis and with
the model error $\delta=8\,000$ and the respective thresholds are 
0.24, 0.4, 0.2 and 0.18.
To reach convergence, we had to let Alg.\,\ref{fig:Condat.algorithm} run for 500 iterations.

\newlength{\imageWidth}
\setlength{\imageWidth}{0.35\textwidth}
\begin{figure*}
    \includegraphics[width=\imageWidth]{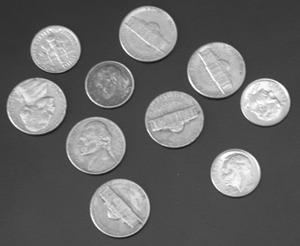}%
    \hfill
    \includegraphics[width=\imageWidth]{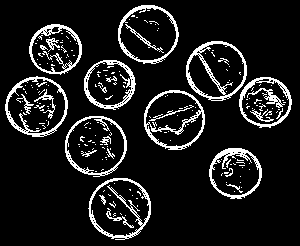}\\[1mm]
    \includegraphics[width=\imageWidth]{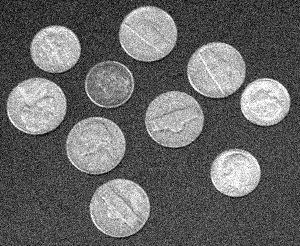}%
    \hfill
    \includegraphics[width=\imageWidth]{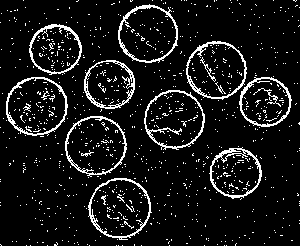}\\[1mm]
    \includegraphics[width=\imageWidth]{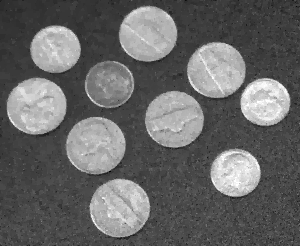}%
    \hfill
    \includegraphics[width=\imageWidth]{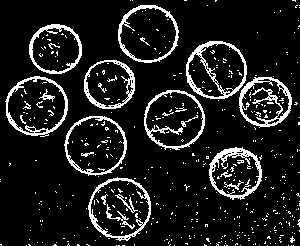}\\[1mm]
    \begin{minipage}[b]{\imageWidth}
        \includegraphics[width=.32\imageWidth]{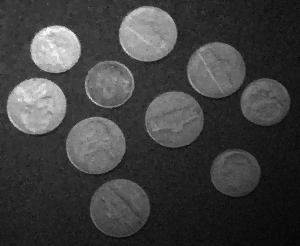}
        \hfill
        \includegraphics[width=.32\imageWidth]{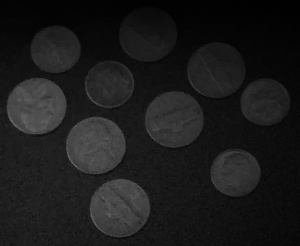}
        \hfill
        \includegraphics[width=.32\imageWidth]{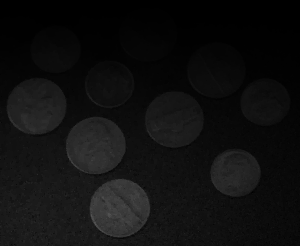}\\[1mm]
        \includegraphics[width=.32\imageWidth]{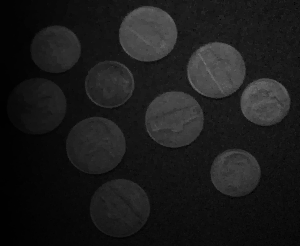}
        \hfill
        \includegraphics[width=.32\imageWidth]{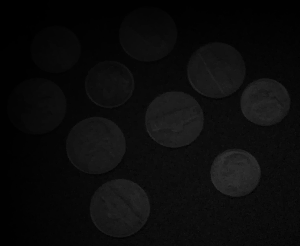}
        \hfill
        \includegraphics[width=.32\imageWidth]{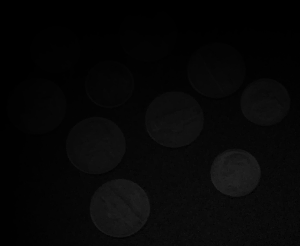}\\[1mm]
        \includegraphics[width=.32\imageWidth]{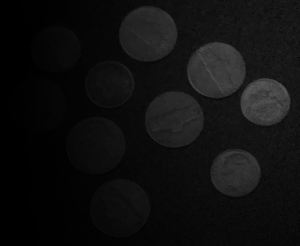}
        \hfill
        \includegraphics[width=.32\imageWidth]{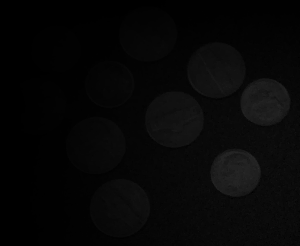}
        \hfill
        \includegraphics[width=.32\imageWidth]{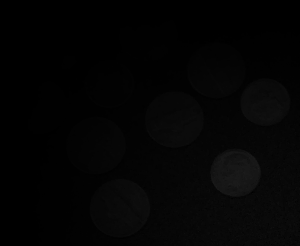}
    \end{minipage}
    \hfill
    \includegraphics[width=\imageWidth]{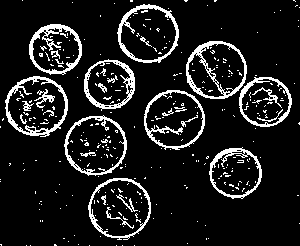}
    \caption{Results of the edge-detecting experiment.
        The left column shows the images, from top to bottom:
        original;
				noisy ($\sigma=20$);
				synthesized by taking $\hat{\y} = \P \hat\x$ (i.e.\ denoised);
				parameter maps $\hat\x_{00}, \hat\x_{01},\ldots, \hat\x_{22}$ (from left to right and from top to bottom, in absolute values and jointly scaled to fit the black--white range).
        The right column shows the respective edge maps obtained by thresholding gradient images.
        Only the last image (which gives best results) is constructed from the parameter maps $\hat\x_{k\ell}$.
    }
    \label{fig:Coins.processed}
\end{figure*}

\begin{figure*}
	\includegraphics[width=\imageWidth]{fig/fig_coins_orig}%
	\hfill
	\includegraphics[width=\imageWidth]{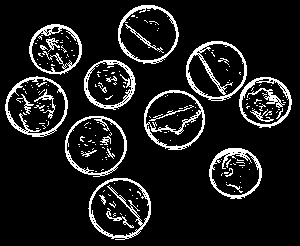}\\[1mm]
	\includegraphics[width=\imageWidth]{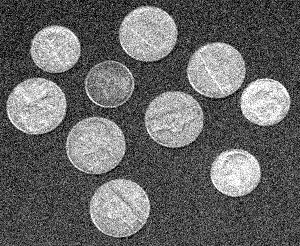}%
	\hfill
	\includegraphics[width=\imageWidth]{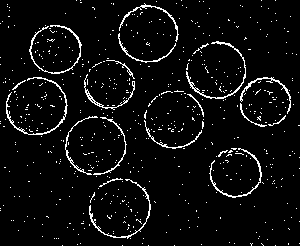}\\[1mm]
	\includegraphics[width=\imageWidth]{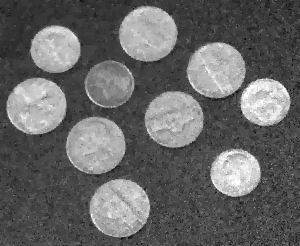}%
	\hfill
	\includegraphics[width=\imageWidth]{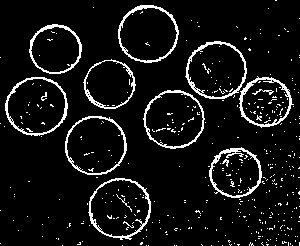}\\[1mm]
	\begin{minipage}[b]{\imageWidth}
		\includegraphics[width=.32\imageWidth]{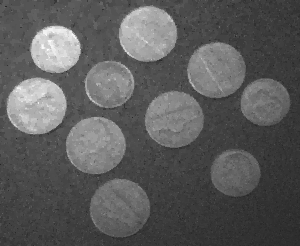}
		\hfill
		\includegraphics[width=.32\imageWidth]{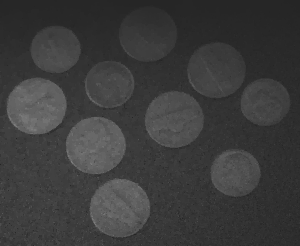}
		\hfill
		\includegraphics[width=.32\imageWidth]{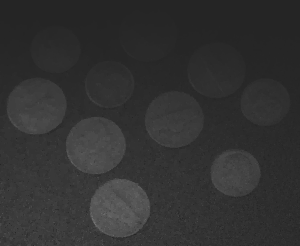}\\[1mm]
		\includegraphics[width=.32\imageWidth]{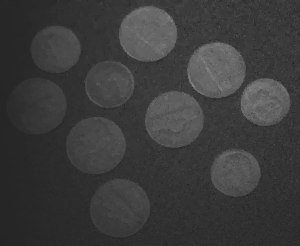}
		\hfill
		\includegraphics[width=.32\imageWidth]{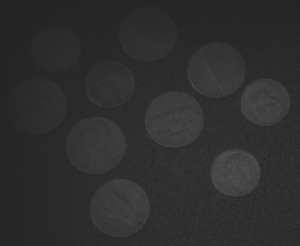}
		\hfill
		\includegraphics[width=.32\imageWidth]{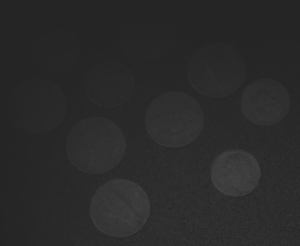}\\[1mm]
		\includegraphics[width=.32\imageWidth]{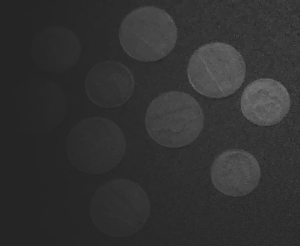}
		\hfill
		\includegraphics[width=.32\imageWidth]{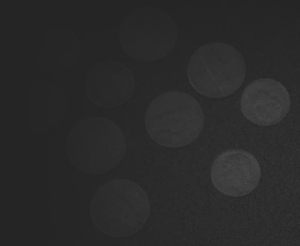}
		\hfill
		\includegraphics[width=.32\imageWidth]{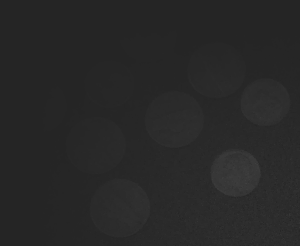}
	\end{minipage}
	\hfill
	\includegraphics[width=\imageWidth]{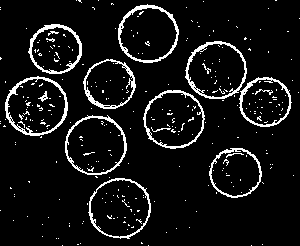}
	\caption{Results of the experiment, similar to the one presented in Fig.\,\ref{fig:Coins.processed}, but for noise standard deviation $\sigma=30$.
	}
	\label{fig:Coins.processed.more.noise}
\end{figure*}

\begin{figure*}
	\includegraphics[width=\imageWidth]{fig/fig_coins_orig}%
	\hfill
	\includegraphics[width=\imageWidth]{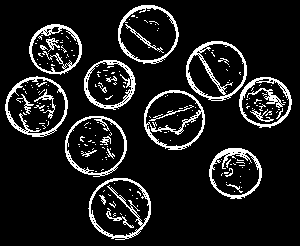}\\[1mm]
	\includegraphics[width=\imageWidth]{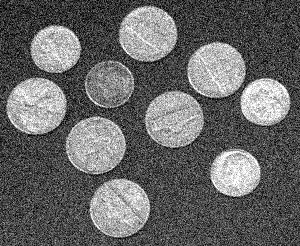}%
	\hfill
	\includegraphics[width=\imageWidth]{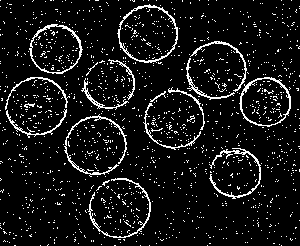}\\[1mm]
	\includegraphics[width=\imageWidth]{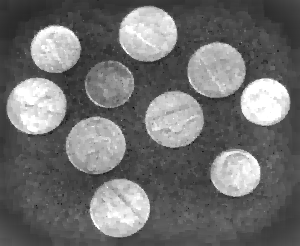}%
	\hfill
	\includegraphics[width=\imageWidth]{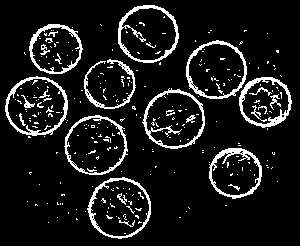}\\[1mm]
	\begin{minipage}[b]{\imageWidth}
		\includegraphics[width=.32\imageWidth]{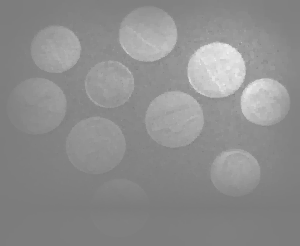}
		\hfill
		\includegraphics[width=.32\imageWidth]{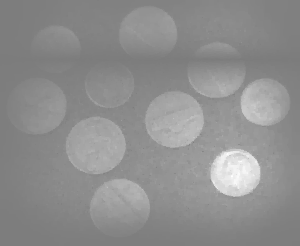}
		\hfill
		\includegraphics[width=.32\imageWidth]{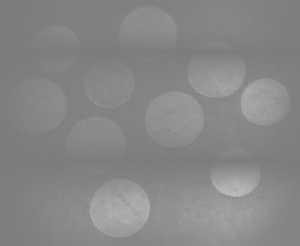}\\[1mm]
		\includegraphics[width=.32\imageWidth]{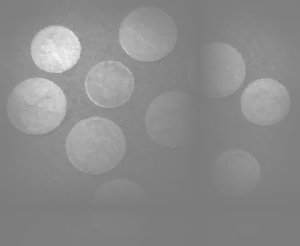}
		\hfill
		\includegraphics[width=.32\imageWidth]{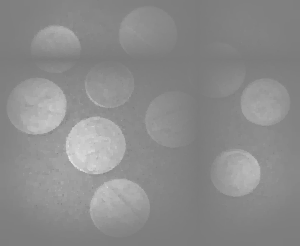}
		\hfill
		\includegraphics[width=.32\imageWidth]{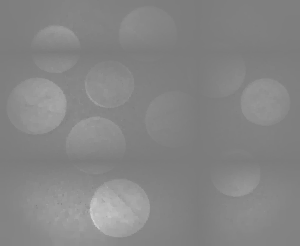}\\[1mm]
		\includegraphics[width=.32\imageWidth]{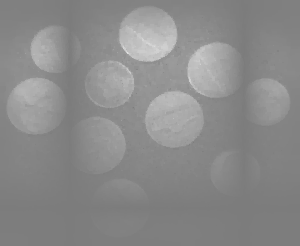}
		\hfill
		\includegraphics[width=.32\imageWidth]{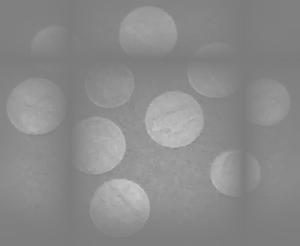}
		\hfill
		\includegraphics[width=.32\imageWidth]{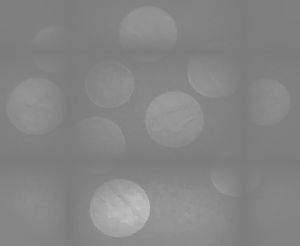}
	\end{minipage}
	\hfill
	\includegraphics[width=\imageWidth]{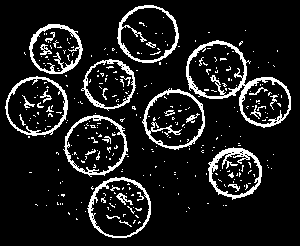}
	\caption{Results of the experiment for $\sigma=30$ as in Fig.\,\ref{fig:Coins.processed.more.noise},
	but now with orthonormal polynomial basis.
	}
	\label{fig:Coins.processed.more.noise.orthonormal}
\end{figure*}

\setlength{\imageWidth}{0.3\textwidth}
\begin{figure*}
	\includegraphics[width=\imageWidth]{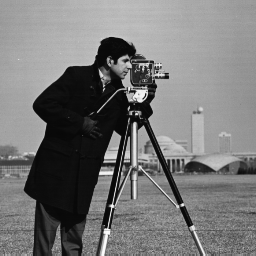}%
	\hfill
	\includegraphics[width=\imageWidth]{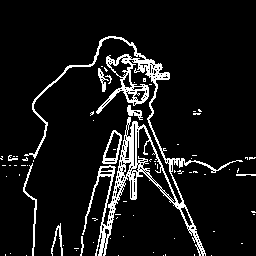}\\[1mm]
	\includegraphics[width=\imageWidth]{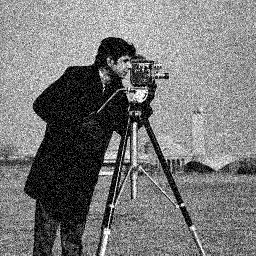}%
	\hfill
	\includegraphics[width=\imageWidth]{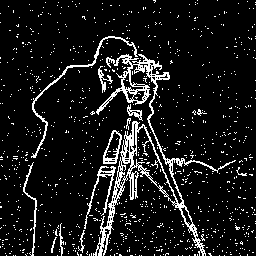}\\[1mm]
	\includegraphics[width=\imageWidth]{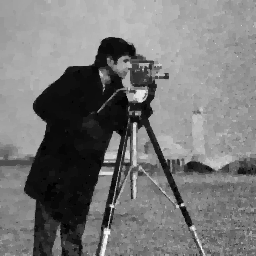}%
	\hfill
	\includegraphics[width=\imageWidth]{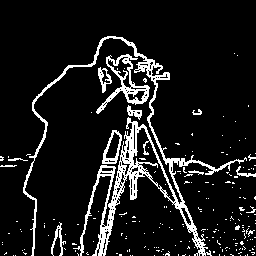}\\[1mm]
	\begin{minipage}[b]{\imageWidth}
		\includegraphics[width=.32\imageWidth]{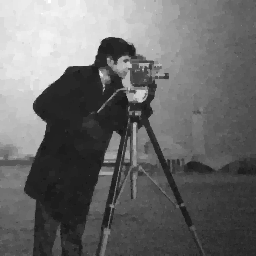}
		\hfill
		\includegraphics[width=.32\imageWidth]{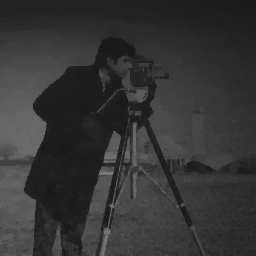}
		\hfill
		\includegraphics[width=.32\imageWidth]{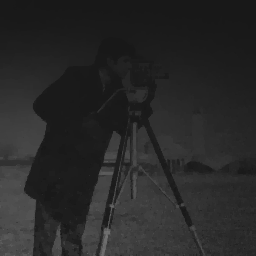}\\[1mm]
		\includegraphics[width=.32\imageWidth]{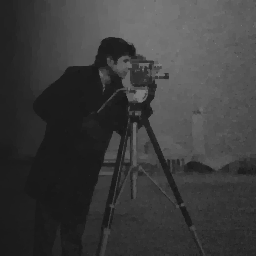}
		\hfill
		\includegraphics[width=.32\imageWidth]{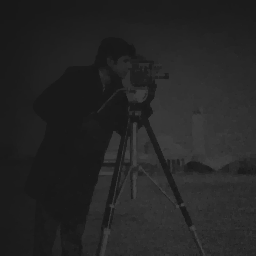}
		\hfill
		\includegraphics[width=.32\imageWidth]{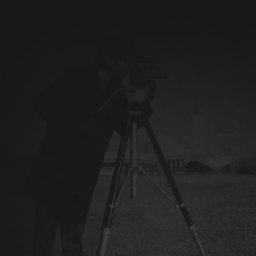}\\[1mm]
		\includegraphics[width=.32\imageWidth]{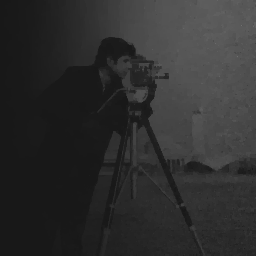}
		\hfill
		\includegraphics[width=.32\imageWidth]{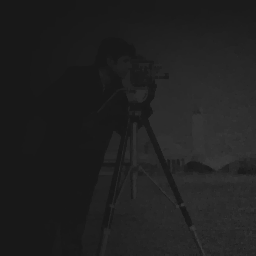}
		\hfill
		\includegraphics[width=.32\imageWidth]{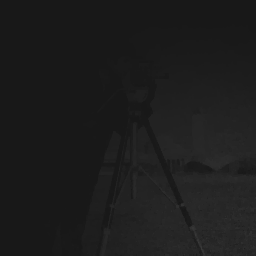}
	\end{minipage}
	\hfill
	\includegraphics[width=\imageWidth]{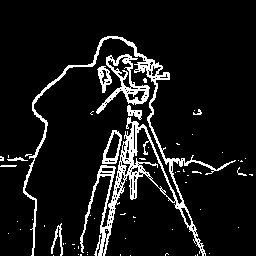}
	\caption{Results of the experiment for the cameraman image, $\sigma=30$.
	}
	\label{fig:Cameraman.processed}
\end{figure*}

\begin{figure*}
	\includegraphics[width=\imageWidth]{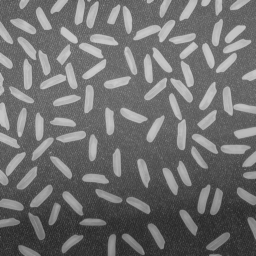}%
	\hfill
	\includegraphics[width=\imageWidth]{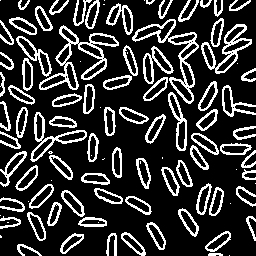}\\[1mm]
	\includegraphics[width=\imageWidth]{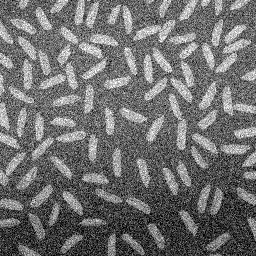}%
	\hfill
	\includegraphics[width=\imageWidth]{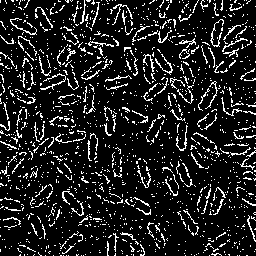}\\[1mm]
	\includegraphics[width=\imageWidth]{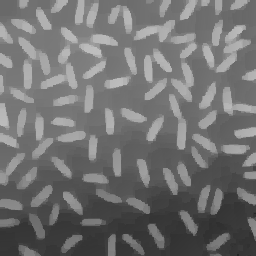}%
	\hfill
	\includegraphics[width=\imageWidth]{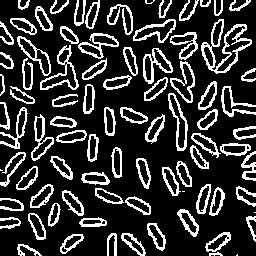}\\[1mm]
	\begin{minipage}[b]{\imageWidth}
		\includegraphics[width=.32\imageWidth]{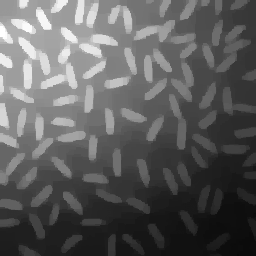}
		\hfill
		\includegraphics[width=.32\imageWidth]{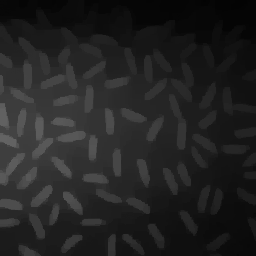}
		\hfill
		\includegraphics[width=.32\imageWidth]{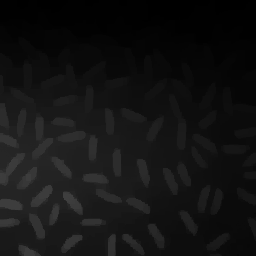}\\[1mm]
		\includegraphics[width=.32\imageWidth]{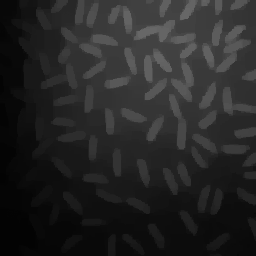}
		\hfill
		\includegraphics[width=.32\imageWidth]{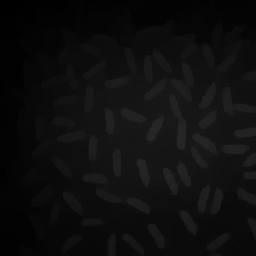}
		\hfill
		\includegraphics[width=.32\imageWidth]{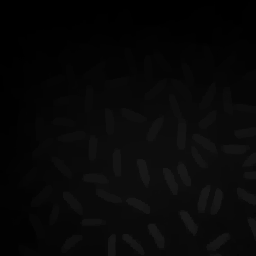}\\[1mm]
		\includegraphics[width=.32\imageWidth]{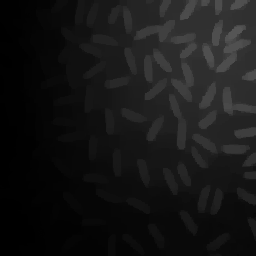}
		\hfill
		\includegraphics[width=.32\imageWidth]{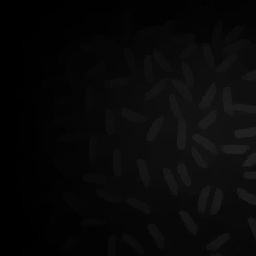}
		\hfill
		\includegraphics[width=.32\imageWidth]{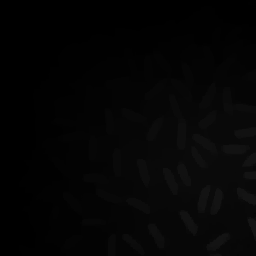}
	\end{minipage}
	\hfill
	\includegraphics[width=\imageWidth]{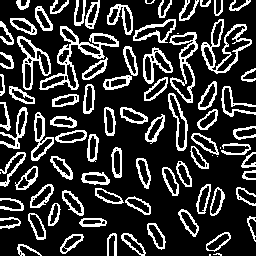}
	\caption{Results of the experiment for the rice image, $\sigma=30$.
	}
	\label{fig:Rice.processed}
\end{figure*}

\subsection{Discussion}
It is clear that with both of our methods we were able to indicate the edges significantly better than the standard detection kernels.
%
Note, however, that the price paid is vastly longer computation times.
The detection based on parameter maps has beaten the simpler synthesizing approach.


\section{Software and data}
The data and the implementation for Matlab can be downloaded from the URL
\url{http://www.utko.feec.vutbr.cz/~rajmic/software/edge_detect.zip}.
The archive is prepared in the sense of reproducible research; note, however,
that the noise is generated each time the script is run,
and that the orthonormal polynomials in experiment \ref{sec:coins.ortho} are based on randomness too.

\section{Conclusion}
We have presented a~model based on overcomplete piecewise-polynomial image assumption with group sparsity.
The proposed model provides the basis for two simple methods for identifying image edges.
We have shown that the two methods resolve image edges robustly with respect to noise,
and that they exhibit a~significant improvement over the performance of the standard Sobel edge detecting kernel.

%

{
\inputencoding{cp1250}
\bibliographystyle{IEEEtran}
\bibliography{literatura}      
}

%
\end{document}